\documentclass[11pt]{amsart}
\usepackage{graphicx,xcolor}
\usepackage{amssymb,amsmath,amsgen,amsopn,amsxtra,amsbsy,amscd}
\usepackage{mathtools}
\usepackage{subcaption}
\numberwithin{equation}{section}

\def\rw{\rightarrow}

\def\({\left(}
\def\){\right)}
\def\nx{\nabla}
\def\Nx{\nabla_x}
\def\Cal{\mathcal}
\def\t{\tilde}
\def \Nx{\nabla}

\def\Om{\Omega}
\def\eb{\varepsilon}

\def\R {\mathbb{R}}

\newcommand{\be}{\begin{equation} }
\newcommand{\ee}{\end{equation} }

\def\LL {{\mathcal L}}

\def \l {\langle}
\def \r {\rangle}
\def \p {\partial}

\def \and{\qquad\text{and}\qquad}

\def\Dx{\Delta_x}
\def\({\left(}
\def\){\right)}
\def\Nx{\nabla}

\def\eb{\varepsilon}
\def\Cal{\mathcal}
\def\eb{\varepsilon}

\def\Om{\Omega}

\def\la{\lambda}

\def\R {\mathbb{R}}

\def\LL {{\mathcal L}}

\def \l {\langle}
\def \r {\rangle}

\def \p {\partial}

\def\Dx{\Delta}

\newtheorem{proposition}{Proposition}[section]
\newtheorem{theorem}[proposition]{Theorem}

\newtheorem{lemma}[proposition]{Lemma}
\theoremstyle{definition}
\newtheorem{definition}[proposition]{Definition}
\newtheorem{remark}[proposition]{Remark}

\numberwithin{equation}{section}

\begin{document}
\title[Chevron Pattern Equations]{Backward behavior and determining  functionals for chevron pattern equations}

\author{ V. K. Kalantarov$^{1}$, H. V. Kalantarova$^{2}$ , and O. Vantzos$^{3}$}
\address{$^{1}$Department of Mathematics, Ko{\c c} University, Rumelifeneri Yolu, Sar{\i}yer, Istanbul, Turkey}
\address{$^{2}$Springer-Verlag, Heidelberg, Germany}
\address{$^{3}$Vantzos Research SMPC, Athens, Greece}


\date{\today}

\maketitle

\begin{abstract}
The paper is devoted to the study of the backward behavior  of solutions of the initial boundary value problem for the chevron pattern equations  under homogeneous Dirichlet's boundary conditions. We prove that, as $t\rightarrow \infty$, the asymptotic behavior of solutions of the considered problem is completely determined by the dynamics of a finite set of functionals. Furthermore, we provide numerical evidence for the blow-up of certain solutions of the backward problem in finite time in 1D.
\end{abstract}

\maketitle

\section{Introduction}

We consider  the system of equations
\begin{eqnarray}
\label{a0}
&&\tau \partial_{t}A=A+\Dx A - \phi^2A-|A|^2A-2i c_1\phi \partial_{y}A+i\beta A\partial_{y}\phi, \\
\label{a1}
&&\partial_{t}\phi=-\mathbb L \phi-h\phi +\phi\lvert A\rvert^{2}-c_{2}\mbox{Im}\left[ A^{\ast}\partial_{y}A\right],\quad (x,y)\in\Omega,\quad t>0,
\end{eqnarray}
where the linear operator $\mathbb L$  is given by
$$
\mathbb L:= -D_1\partial^{2}_{x}-D_2\partial^{2}_{y},
$$
and $\Omega\subset\mathbb{R}^{2}$ is a bounded domain with a sufficiently smooth boundary $\partial\Omega$. The system is complemented by the following boundary conditions
\be\label{a2}
A\Big|_{\p\Om}=0,\quad \phi\Big\vert_{\partial\Omega}=0,\quad t>0,
\ee
and the initial conditions
\be\label{a3}
A\Big|_{t=0}=A_{0},\quad \phi\Big\vert_{t=0}=\phi_{0},\quad (x,y)\in\Omega.
\ee
In the context of this study,
\be\label{cd1}D_1>0,\ D_2>0,\ c_1\geq0,\ c_2\geq0,\ h\geq0,\ \beta \in \R
\ee 
are given parameters, the complex valued function $A$, where $A^{\ast}$ denotes its complex conjugate, and the real valued function $\phi$ are unknown functions. This coupled system of equations was proposed by Rossberg et al. \cite{Rossberg_diss}, \cite{Rossberg1996weakly} to model the formation and evolution of patterns in nematic liquid crystals subject to convection driven by an electric field, sometimes accompanied with magnetic field. These works were motivated by a number of experimental studies \cite{Kai1996}, \cite{Richter1995}, \cite{Nasuno1992}.
  \begin{figure}[h]
  \centering
  \includegraphics[width=0.7\textwidth]{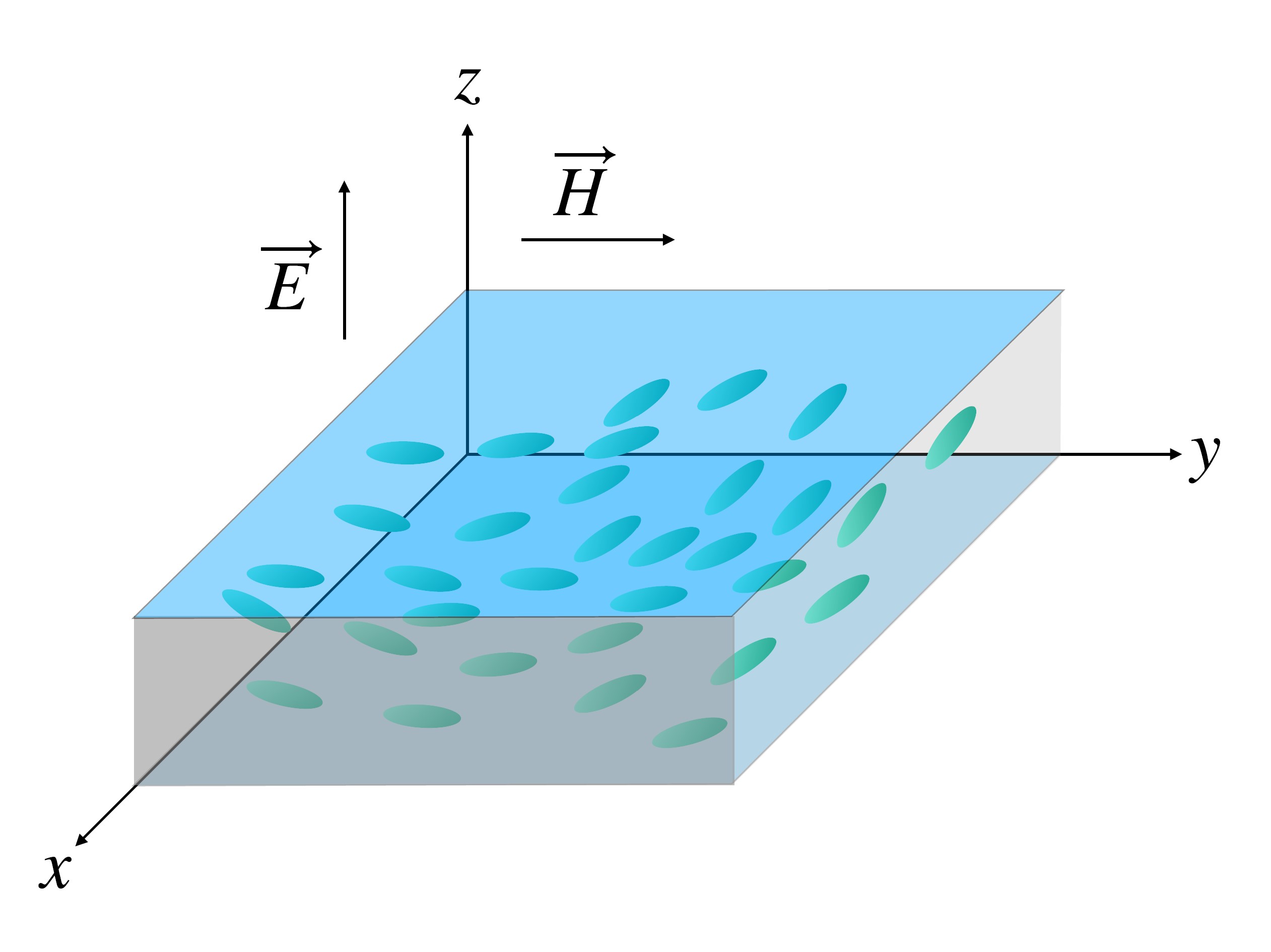} 
  \caption{\footnotesize{Cross-section of a nematic liquid crystal contained between two parallel transparent plates. Nematic liquid crystals are composed of rod-like molecules, that are free to flow and orient themselves in 3D space, when external electric $\vec{E}$ and magnetic fields $\vec{H}$ are applied.}}
  \label{fig:cross_section}
 \end{figure}
A typical experimental setup used in these investigations consists of a nematic liquid crystal confined between two planar transparent glass plates, placed parallel at a distance of approximately $50 \mu m$ from each other. Then AC-voltage of varying frequency is applied across the plates, sometimes in conjunction with a magnetic field parallel to the plates. The results of the described paradigmatic experiment, revealed that through use of either microscopes (with or without polarizers) or via laser diffraction various pattern formations are observable: Williams domain, Grid-Pattern, chevron patterns, etc.. The formation of these patterns is governed by experimentally easily controllable parameters, the amplitude and the frequency of the applied voltage and the magnitude of magnetic field. The advantage of liquid crystals over other fluids in these experiments besides its transparency and its phase transition temperatures being close to the room temperature, which makes detection of patterns easier, \cite{Buka2001}, is the fact that their molecules are locally oriented along desired direction. In \cite{Rossberg1998}, it is reported that the described system has a tendency to form chevron structures after a transition to defect chaos and that chevron patterns are generally observed when the isotropy of the system is spontaneously broken, this effect is achieved either by increasing the voltage or by applying a weak magnetic field. 

The model \eqref{a0}-\eqref{a1} is obtained by coupling the real Ginzburg-Landau equation, where $A$ represents the phase, direction and amplitude of the periodic pattern and equation \eqref{a1}, where $\phi$ represents the angle between director and the x-axis and also by using symmetries, smallness and regularity observations from experiments as assumptions. In this model the parameter $\tau$ is a function of the various physical timescales of the problem, $\beta$ measures the interaction between the gradient of the director field and the phase of the rolls, the parameters $D_{1}$ and $D_{2}$ are the coefficients of the anisotropic diffusion of the director field for the liquid crystal particles, the dampening parameter $h$ measures the tendency of the director field to align with the magnetic field $H$, corresponding to $\phi=0$ and the parameters $c_{1}$ and $c_{2}$ control the torque that the director field and the wave vector of the rolls exert on each other.

In \cite{Rossberg1996weakly}, Rossberg et al. generated chevron patterns in simulations of \eqref{a0}-\eqref{a1}. For more numerical studies of the model refer to \cite{Sakaguchi2009}, \cite{Komineas2003}.

The mathematical model \eqref{a0}-\eqref{a1} has also been studied qualitatively. In \cite{KKV}, the authors prove the following well-posedness result:

 \begin{theorem}\label{exun} (\cite{KKV})
If $c_1\in [0,1)$ or $c_1\ge 2c_2$, $h>0$, then the  initial boundary value problem \eqref{a0}-\eqref{a3}, where $A_{0},\phi_{0}\in L^{2}(\Omega)$, has a  unique  weak solution
\begin{equation}
A,\phi \in C([0,T]; L^2(\Om)) ; H_0^1(\Om)),\quad \forall T>0,  \ \
\end{equation}
such that
\be\label{aest1}
\|\phi(t)\|\le M_0, \ \ \|A(t)\|\le M_0, \ \ \forall t>0,
\ee
 and
\be\label{aest2}
\int_0^T\|\nx \phi(t)\|^{2}dt \le M_T, \ \  \int_0^T\|\nx A(t)\|^{2}dt\le M_T, \ \forall T>0,
\ee
where $M_{0}$ denotes a generic constant depending only on $\lVert A_{0}\rVert$, $\lVert\phi_{0}\rVert$ and $\lvert\Omega\rvert$, and where $M_{T}$ is a generic constant, which depends also on T, besides $\lVert A_{0}\rVert$, $\lVert\phi_{0}\rVert$ and $\lvert\Omega\rvert$. In other words this problem generates a continuous semigroup $S(t$), $t\ge0,$ in the phase space $V^{0}:= L^2(\Omega)\times L^2(\Omega).$
Moreover, this semigroup is bounded dissipative in the phase space $V^{0}$.
\end{theorem}

\noindent In \cite{KKV}, the authors also proved that the problem \eqref{a0}-\eqref{a3} possesses a global attractor, i.e., there exists 
a compact, invariant set that attracts uniformly each bounded set of the phase space $V^{0}$. Furthermore, it is proved that this semigroup has a global attractor that belongs to $V^{1}:=H^{1}_{0}(\Omega)\times H^{1}_{0}(\Omega)$. Later, this result is improved in \cite{KKV2021}, in which it is proved that under the same assumptions of the Theorem~\ref{exun} the semigroup\\  $S(t):V^{0}\rightarrow V^{0}$, $t\ge0,$ generated by the problem \eqref{a0}-\eqref{a3} possesses an exponential attractor.

In \cite{KKV2021}, the authors also prove the global stabilization of the zero steady state of the system \eqref{a0}-\eqref{a1} by finitely many Fourier modes:

\begin{theorem}
The solution $[A(t),\phi(t)]$  of the problem 
\begin{eqnarray}
\label{chp1}&&\tau\p_t A=A+\Delta A -\phi^2A-|A|^2A-2i c_1\phi\partial_y A
+i\beta A\partial_y \phi,\\
&&\quad\quad\quad\quad\quad\quad\quad\quad\quad\quad\quad\quad\quad\quad\quad\quad-\mu\sum\limits_{k=1}^N(A,w_k)w_k\nonumber\\
\label{chp2}&&\partial_t \phi=-\mathbb L \phi -h\phi +\phi\lvert A\rvert^{2}-c_{2}\mbox{Im}
\left[ A^{\ast}\partial_y A\right],\\
\label{chp3}&&A\big|_{\p \Om}=\phi\big|_{\p \Om}=0, \ \ A\big|_{t=0}=A_0, \ \phi\big|_{t=0}=\phi_0,
\end{eqnarray}
tends to zero stationary state  with an exponential rate in $V^{1}$, if 
\begin{equation}
\mu\geq 1\quad\mbox{ and }\quad\lambda_{N+1}^{-1}<\delta_{0},\nonumber
\end{equation}
where $0<\lambda_{1}\leq\lambda_{2}\leq\ldots\leq\lambda_{N+1}\leq\ldots$ are the eigenvalues of the Laplace operator under the homogeneous Dirichlet boundary condition,
\begin{equation}
\delta_{0}:=2(1-c_{1})/(2+c_{2})\quad\mbox{and}\quad c_{1}<1.\nonumber
\end{equation}
\end{theorem}

In this work, our primary objective is to demonstrate that the  long time behavior of solutions to the initial boundary value problem for the system of equations \eqref{a0}-\eqref{a1} can be uniquely determined by finite number of determining functionals. Furthermore, we study the backward in time behavior of solutions to the system \eqref{a0}-\eqref{a1}.

After covering some preliminaries in section \ref{sec:preliminaries}, we devote section \ref{sec:backward} of the paper to the backward in time behavior of solutions to the system \eqref{a0}-\eqref{a1}. There are numerous  papers dedicated to the backward behavior of solutions to nonlinear evolution equations (see, e.g.
 Bardos and Tartar \cite{BaTa}, Dascaliuc \cite{Da}, Kukavica and Malcok \cite{KuMal}, Guo and Titi \cite{GuTi}, Vukadinovic \cite{Vu}  and references therein). In this section, we establish the backward uniqueness of the problem \eqref{a0}-\eqref{a3} and show that some solutions of the backward problem blow up in a finite time.
We follow this with further numerical evidence of finite-time blow up in  section \ref{sec:numerics}, using a convex-concave splitting scheme with explicit tracking of blow up times during each discrete time step.

Section \ref{sec:functionals} of the paper is dedicated to studying the existence of determining functionals for the system of chevron pattern equations \eqref{a0}-\eqref{a1}. The intensive study over the last decades of finite  dimensional asymptotic behavior of dissipative dynamical systems generated by nonlinear PDEs has been inspired by the pioneering works of Foias and Prodi \cite{FP} and Ladyzhenskaya \cite{LA1}. They proved that asymptotic behavior of solutions to the initial boundary value problem for 2D Navier-Stokes equations is completely determined by asymptotic behavior of first $N$ Fourier modes (for $N$ large enough). For  subsequent investigations of  the finite-dimensional behavior of  dynamical systems generated by the initial boundary value problem for dissipative  PDEs (existence of finite-dimensional global attractors,  inertial manifolds, exponential attractors, finite number of  modes, nodes and local volume elements determining asymptotic behavior of solutions) we refer to \cite{BV}, \cite{FMRT}, \cite{HaRo}, \cite{Te} and the references therein.
The concept of determining functionals (determining interpolant operators) is introduced in \cite{CJT1}, \cite{CJT2}. This concept provided a unified approach for investigating parameters that uniquely determine the asymptotic behavior of solutions to dissipative PDEs. A further extension  of this unified approach for studying  the  long-time behavior of various  nonlinear dissipative PDEs was developed  in subsequent investigations of the long-time behaviour of  solutions to initial boundary value problems for  Navier-Stokes equations, systems of nonlinear parabolic equations, weakly damped nonlinear wave equations, nonlinear wave equation with nonlinear damping,  von Karman equations   and some other dissipative systems (see, e.g., Chueshov \cite{Chu1}, \cite{Chu2}, \cite{Chu3},  Chueshov and Kalantarov \cite{ChuKal}, Chueshov and Lasiecka \cite{ChuLas},
 Holst and Titi \cite{HT} and references therein.)
 
The main result of this section is Theorem \ref{DF2}, which establishes the existence of a finite number of functionals that uniquely determine the asymptotic behavior of solutions to the problem \eqref{a0}-\eqref{a3}.\\

\section{Preliminaries}\label{sec:preliminaries}

In this section, we introduce the notation that will be used in what follows and recapitulate some definitions and known results.

\begin{itemize}
\item $C$ will denote a generic positive constant,
\item $(\cdot,\cdot):=(\cdot,\cdot)_{L^2(\Om)}$, $\| \cdot\|:=\| \cdot\|_{L^2(\Om)}$,
\item $V^{0}:=L^2(\Omega) \times L^2(\Omega )$ is the Hilbert
space with the inner  product
$$\l [u_1,v_1], [u_2,v_2]\r_{V^0}:=( 
u_1, u_2)+(  v_1, v_2),$$
\item $V^1:=H^1_{0}\times H^1_{0}$ denotes the Hilbert space with  the inner product
$$\l [u_1,v_1], [u_2,v_2]\r_{V^{1}}:=( \Nx
u_1,\Nx u_2)+( \nabla v_1,\nabla v_2).
$$
\end{itemize}

\begin{definition}\label{D1}( \cite{CH} ) Let $V$ and $H$ be reflexive Banach
spaces and $V$ be continuously and densely embedded into $H$. The
{\it completeness defect}  of a set ${\Cal L}=\{\ell_j:j=1,\ldots,N\}$ of
linear functionals on $V$ with respect to $H$ is the following number
$$\epsilon_{{\Cal L}}(V,H):=\sup \{\|
w\|_H: w\in V,\,\ell_j(w)=0,\,\ell_j\in {\Cal L},\,\| w\|_V\leq 1\}.
$$
$\epsilon_{{\Cal L}}(V,H)=0$ if and only if the class of functionals
${\Cal L}$ is complete in $V$, i.e., the condition $\ell_j(w)=0$ for
all $\ell_j\in {\Cal L}$ implies $w=0$.
\end{definition}

\begin{theorem}\label{Ch}  ( \cite{CH} ) Let $\epsilon_{\Cal L}= \epsilon_{\Cal
L}(V,H)$ be the completeness defect of a set ${\Cal L=\{\ell_{j}:j=1,\ldots,N\}}$ of linear
functionals on $V$ with respect to $H$. Then there exists a constant
$C_{\Cal L}$ such that
$$
\| w\|_H\leq C_{\LL}\eta_{\LL}(w)+ \epsilon_{\Cal
L}\,\| w\|_V \qquad \mbox{for any} \quad w\in V,$$
where $\eta_{\LL}(w)=\max\{|\ell_j(w)|: j=1,\ldots,N\}$.
\end{theorem}

\begin{lemma}\label{uGr} (see, e.g., \cite{Te}) Let $g,h$ and $y$ be three positive locally integrable functions on the interval $(t_0,\infty)$, such that $y'$ (or equivalently $\frac{dy}{dt}$) is locally integrable on the same interval. Suppose that the following inequalities are satisfied
$$
y'(t)\le g(t)y(t)+h(t), \ \ t\in [t_0,\infty),
$$
$$
\int_t^{t+r}g(s)ds\le a_1, \  \int_t^{t+r}h(s)ds\le a_2, \ \int_t^{t+r}y(s)ds\le a_3,  \ \ \forall  t\in [t_0,\infty),
 $$
where $r,a_1,a_2,a_3$ are positive constants.\\
Then
 $$
 y(t+r)\le (\frac{a_3}r+a_2)e^{a_1}, \ \ t\in [t_0,\infty).
 $$
\end{lemma}

\begin{theorem}\label{tabs} \mbox{(see, e.g.~\cite{Te})} Suppose that
$\mathbf A$  is a positive definite, unbounded, self adjoint linear operator in H with domain $V^{2}: = D(\bf{A})$ and let the nonlinear operator $\mathbf{B}(\cdot):  (0,T)\times V^{1}\rw H , $ with $V^{1}:=D(\mathbf{A}^{\frac12})$ satisfy the condition
\be\label{abs1}
\|\mathbf{B}(\mathbf{u}(t))\|_H\le k(t)\|\mathbf{u}(t)\|_{V^{1}},\quad k(t)\in L^2(0,T).
\ee 

If
\be\label{abs2}
\mathbf{u}\in L^\infty(0,T; V^{1})\cap L^2(0,T; V^{2})
\ee
is a solution of the equation 
\be\label{abs}
\mathbf{u}'(t)+\mathbf{A}\mathbf{u}(t)= \mathbf{B}(\mathbf{u}(t)), \ t\in (0,T),
\ee
that satisfies the condition
 $\mathbf{u}(T)=0$, then $\mathbf{u}(t)=0, \ \forall t\in [0,T].$
\end{theorem}
 
\begin{definition} A set of continuous linear functionals $\mathcal{L}=\{l_{j}:j=1,\ldots,N\}$ on $V$ is said to be a set of asymptotically $(V,H,\mathcal{W})-$determining functionals for the problem \eqref{abs} if for any two solutions $u,v\in\mathcal{W}$ the condition
\begin{equation*}
\lim_{t\rightarrow \infty}\int_{t}^{t+1}\lvert l_{j}(u(\tau))-l_{j}(v(\tau))\rvert^{2}d\tau=0\mbox{ for }j=1,\ldots,N,
\end{equation*}
implies that
\begin{equation*}
\lim_{t\rightarrow\infty}\lVert u(t)-v(t)\rVert_{H}=0.
\end{equation*}
\end{definition}

\begin{theorem}\label{Det1} (\cite{Chu1}, Theorem 4.1.) Suppose that \begin{itemize}
\item $\mathbf A$ is a positive definite operator  with a discrete spectrum in a Hilbert space $H$, the basis $\{w_k\}$ of which consists of eigenvectors of the operator $\mathbf{A}$:
$$
\mathbf Aw_k=\la_k w_k, \ \ (w_k,  w_j)=\delta_{kj}, \ 0<\la_1\le\la_2\le \cdots, \ \ \lim_{k\rightarrow\infty} \la_{k}=\infty,
$$
\item the nonlinear operator $\mathbf B(\cdot): D(\mathbf A^{\frac12})\rw D(\mathbf A^{-\gamma}), \ 0\le \gamma <\frac12$  satisfies the condition
\be\label{Lip1}
\|\mathbf A^{-\gamma}(\mathbf B({\mathbf u}_1)-\mathbf B({\mathbf u}_2))\|\le M(\rho)\|\mathbf A^{\frac12}(\mathbf u_1-\mathbf u_2)\|, 
\ee
for all ${\mathbf u}_j\in \mathbf D( {\mathbf A}^{\frac12})$ such that $\|\mathbf A^{\frac12}\mathbf{u_j}\|\le \rho$, $\rho>0,$ and $M(\rho)>0$ is a continuous function on $\R^+,$\\
\item the problem 
\be\label{abs}
\begin{cases}
\mathbf{u}'(t)+\mathbf{A}\mathbf{u}(t)= \mathbf{B}(\mathbf{u}(t)),  \ \ t>0,\\
 \mathbf u(0)=\mathbf u_0, 
\end{cases}
\ee
 has a unique solution in the class of functions $\mathcal W=C([0,\infty);H)\cap  C([0,\infty); D(A^{\frac12}))$ and it is point dissipative , i.e., 
$$
\|\mathbf A^{\frac12}\mathbf{u}(t)\|\le R, \ \ \forall t\ge t_0(u),\quad \forall \mathbf{u}(t)\in\mathcal{W}.
$$
\end{itemize}
Suppose also that $\LL=\{\ell_j, \ j=1,\ldots, N\}$ is a set of linear continuous functionals defined on $V^1=D(\mathbf A^\frac12)$. Then $\LL$ is a set of asymptotically  
$(V^1,H,\mathcal W)$- determining functionals of the problem \eqref{abs} whenever
$$
\eb_{\LL }(V^1,H)<\eb_0(\gamma, R)\equiv \sqrt{\frac{1+2\gamma}{1-2\gamma}}\left[(1+2\gamma)M(R)\right]^{\frac1{2\gamma-1}}.
$$
\end{theorem} 

\begin{definition}\label{soln}We call a pair of functions $[A,\phi]$ a weak solution of the initial boundary value problem for the system of equations \eqref{a0}-\eqref{a1} under the homogeneous Dirichlet\rq{}s boundary condition, if  $[A,\phi]\in \mathcal{W}= C([0,\infty);V^{0})\cap C((0,\infty);V^{1})$ and the system of equations \eqref{a0}-\eqref{a1} is satisfied in the sense of distributions.
\end{definition}

\section{Backward behavior}\label{sec:backward}

In this section, we will prove the backward uniqueness of a solution to the problem \eqref{a0}-\eqref{a3} and demonstrate that some solutions of the backward problem blow up in finite time. For 1D chevron pattern equations, we will illustrate that solutions not belonging to the global attractor of the semigroup generated by the initial boundary value problem blow up in a finite time $T>0$.

The backward uniqueness of solution to the problem follows from the corresponding result for an abstract differential equation in a Hilbert space $H$ of the form \eqref{abs}.

Suppose that $[A, \phi]$ and $[\t A, \t \phi]$ are solutions of the system \eqref{a0}-\eqref{a1} corresponding to initial data $[A_0,\phi_0]$ and $[\t A_{0}, \t \phi_{0}]$ respectively, then \\
$[a, \Phi]:=[A-\t A, \phi-\t\phi]$ is a solution of the following system 
\begin{multline}\label{unc1}
\tau \p_t a =a+\Dx a-\phi^2a-\t A(\phi+\t \phi) \Phi
 -|A|^2a-(A a^*+\t A^{\ast} a)\t A\\-
2ic_1[\Phi \p_yA+\t \phi \p_y a]+i\beta[a\p_y\phi+\t A\p_y\Phi],
\end{multline}
\begin{multline}\label{unc2}
\p_t \Phi =-\mathbb{L} \Phi -h\Phi+\Phi|A|^2+\tilde{\phi}(Aa^*+\t A^{\ast}a) -c_{2}\mbox{Im} \left[ a^{\ast}\partial_{y}A+\tilde{A}^{\ast}\partial_{y} a\right].
\end{multline}

We can write this system in the form \eqref{abs} in a Hilbert space $V^{0}=L^2(\Om)\times L^2(\Om)$ with  $\mathbf{u}:=\left[
                   \begin{array}{c}
                      a \\
                     \Phi\\
                   \end{array}
                 \right], $
\be\label{wAh}
 \mathbf{A}:=
                 \left[
    \begin{array}{cc}
      -\frac{1}{\tau}\Delta & 0\\
      0 & \mathbb L \\
    \end{array}
  \right],\quad \mathbf {B}(a,\Phi) :=\left[
                   \begin{array}{c}
                    \frac{1}{\tau} f_1(a,\Phi) \\
                     f_2(a,\Phi)\\
                   \end{array}
                 \right],
\ee
where
\begin{multline}\label{f11}
f_1(a,\Phi):=a-\phi^2a-(\phi^2-\t \phi^2)\t A -|A|^2A+|\t A|^2\t A\\-
2ic_1[\Phi \p_yA+\t \phi \p_y a]+i\beta[a\p_y\phi+\t A\p_y\Phi],
\end{multline}

\be\label{f22}
f_2(a,\Phi):=-h\Phi+\Phi|A|^2+\tilde{\phi}(|A|^2-|\t A|^2)\\ -c_{2}\mbox{Im} \left[ a^{\ast}\partial_{y}A+\tilde{A}^{\ast}\partial_{y} a\right]
\ee

$$
V^1:=D(\mathbf{A}^{\frac12})=H_0^1(\Om)\times  H_0^1(\Om)
$$
and 
\be\label{V1}
V^{2}:=D(\mathbf{A})=H^2(\Om)\cap H_0^1(\Om)\times H^2(\Om)\cap H_0^1(\Om).
\ee
It is clear that the operator $\mathbf{A}$ is a positive definite self-adjoint operator in $V^{0}:=L^2(\Om)\times L^2(\Om)$, with the domain $D(\mathbf{A})$   and the nonlinear term $\mathbf{B}(a,\Phi)$ defined in \eqref{wAh} satisfies condition \eqref{abs1}. Thus the statement of the Theorem \ref{tabs} holds true also for the system \eqref{a0}-\eqref{a1}. In other words, solution of the initial boundary value problem for the  backward chevron pattern equations, i.e.,  solution of the problem
\be\label{bchp1}
- \frac{\partial A}{\partial t}=A+\Dx A - \phi^2A-|A|^2A-2i c_1\phi \frac{\partial A}{\partial y}+i\beta A\frac{\partial \phi}{\partial y},
\ee
\be\label{bchp2}
-\frac{\partial \phi}{\partial t}=D_1\frac{\partial^{2}\phi}{\partial x^{2}}+D_2\frac{\partial^{2}\phi}{\partial y^{2}} -h\phi +\phi\lvert A\rvert^{2}-c_{2}\mbox{Im}\left[ A^{\ast}\frac{\partial A}{\partial y}\right],
\ee
\be\label{bchp3}
 A\big|_{\p \Om}=\phi\big|_{\p \Om}=0, \ \
A\big|_{t=0}=A_0,\ \phi\big|_{t=0}=\phi_0
\ee
is unique.\\
Theorem \ref{tabs} also guarantees backward uniqueness of 1D version of the system \eqref{a0}-\eqref{a1}, i.e.,  uniqueness of a weak solution to the problem
\begin{eqnarray}
\label{A0}&&\p_t A+A+\partial_x^2A-\phi^2A =|A|^2A,\quad x\in (0,L),\quad t>0,\\
\label{A1}&&\phi_t +D_1\partial^2_x\phi-h\phi+|A|^2\phi=0,\quad x\in (0,L),\quad t>0,\\
\label{sec3:A2}&&A\big|_{x=0}= A\big|_{x=L}=\phi\big|_{x=0}=\phi\big|_{x=L}=0,\\ 
\label{sec3:A3}&&A\big|_{t=0}=A_0,\quad \phi\big|_{t=0}=\phi_0.
\end{eqnarray}
It is not difficult to see that  some solutions of the  initial boundary value problem \eqref{A0}-\eqref{sec3:A3}
blow up in finite time.

In fact suppose that $[A,\phi]$ is a local solution of the problem corresponding to the initial data $[A_0,\phi_0]$.
Multiplying the equation \eqref{A0} by $A$ in $L^2$ sense we get
$$
\frac12 \frac{d}{dt}\|A(t)\|^2=- \|A(t)\|^2+ \|\p_xA(t)\|^2+(\phi^2(t),|A(t)|^{2})+\int_0^{L}|A(x,t)|^4dx. 
$$
Thanks to the Cauchy-Schwarz inequality we have
$$
\int_0^L|A(x,t)|^2dx\le \sqrt L \left(\int_0^{L}|A(x,t)|^4dx\right)^{\frac12}.
$$
Thus employing the last inequality and the Poincare inequality we  obtain that the function $\Psi(t)=\|A(t)\|^2$  satisfies the differential inequality
$$
\Psi'(t)\ge -2 \Psi(t) +2\la_1 \Psi(t)+2L^{-1} \Psi^2(t).
$$
It is clear that if $\Psi_0:=\Psi(0)=\|A_0\|^2>(1-\la_1)L$, then the solution of this inequality tends to infinity in a finite time
$$
T_0\le \int_{\Psi_0}^\infty\frac {ds}{2(L^{-1}s^2+(\la_1-1)s)}, $$ i.e., the solution of the problem \eqref{A0} blows up in a finite time.\\

\section{Numerical Blow-up of the Backward System}\label{sec:numerics}

To illustrate numerically the blow-up of $\lVert A\rVert$ in the 1D backward system \eqref{A0}, we consider the case where $\phi=0$ 
and note immediately that the resulting equation $\partial_t A = -\partial^2_x A +(\lvert A\rvert^2-1)A$ can not be integrated numerically 
in a straightforward manner because of the 'wrong' sign of the Laplacian $\partial^2_x A$. 
It is a basic result in numerical analysis that a direct numerical integration of this equation would be dominated by spurious high frequency modes.

To progress, we notice that the equation, in its weak form, can be thought of as the \emph{maximizing $L^2$-gradient flow}
\begin{align*}
  &(\partial_t A,u) = \mathcal{F}'(A)(u), \quad u\in H^1_0([0,L]),\\
  &\mathcal{F}(A)=\frac{1}{2}\lVert\partial_x A\rVert^2 + \frac{1}{4}\lVert \lvert A\rvert^2-1\rVert^2 \,.
\end{align*}
Given that the main challenge with integrating this gradient flow numerically is the uncontrollable divergence of the gradient $\partial_x A$,
we consider instead the maximizing gradient flow of the functional $\mathcal{F}$ in the $H^1_\epsilon$-norm 
$$\lVert A\rVert_{H^1_\epsilon}^2 := \lVert A\rVert^2 + \epsilon\lVert\partial_x A\rVert^2,\quad\epsilon>0,$$
which can be written as:
\begin{equation*}
  (\partial_t A,u) + \epsilon(\partial_x\partial_t A, \partial_x u) = \mathcal{F}'(A)(u) \,.
\end{equation*}
This regularization controls the `rate of ascent' of $A$ in $H^1$ instead of $L^2$, 
and so penalizes the rapid divergence of the gradient that accompanies the introduction of spurious oscillatory modes.
In strong form, it yields the regularized PDE
\begin{equation}\label{Aeps}
  \partial_t A - \epsilon \partial_x^2\partial_t A = -\partial^2_x A +(\lvert A\rvert^2-1)A \,.
\end{equation}

Proceeding with the numerical integration of \eqref{Aeps}, we assume that $A(t^k)=A^k$ and apply \emph{operator splitting} to do the integration in two steps:
first integrate the non-linear part
\begin{equation}\label{opsplit1}
  \partial_t \hat A - (\lvert \hat A\rvert^2-1)\hat A = 0,\quad\hat A(t^k)=A^k,
\end{equation}
to calculate $\hat A^{k+1} = \hat A(t^{k+1})$, $t^{k+1}=t^k+\tau$, and then the linear part
\begin{equation}\label{opsplit2}
  (id-\epsilon\partial_x^2)\partial_t A = -\partial_x^2 A,\quad A(t^k)=\hat A^{k+1},
\end{equation}
to update $A$ to the next time step.

The first step \eqref{opsplit1} can be integrated exactly:
$$ \hat A(x,t) = \frac{A^k(x)}{\sqrt{\lvert A^k(x)\rvert^2 + (1-\lvert A^k(x)\rvert^2)e^{2(t-t^k)}}}\,.$$
This blows up for $\lvert A^k(x)\rvert>1$ at the time $$t = t^k + \frac{1}{2}\log\left(\frac{\lvert A^k(x)\rvert^2}{\lvert A^k(x)\rvert^2-1}\right)\,.$$
It follows that
\begin{equation}
  \hat A^{k+1} = \frac{A^k}{\sqrt{\lvert A^k\rvert^2 + (1-\lvert A^k\rvert^2)e^{2\tau}}} \,,
\end{equation}
as long as the (global) time step $\tau$ is chosen so that 
$$\tau<\tau_{blow}=\inf_x\frac{1}{2}\log\left(\frac{\max(1,\lvert A^k(x)\rvert^2)}{\max(1,\lvert A^k(x)\rvert^2)-1}\right)\,.$$
Having selected a suitable $\tau$, in particular ensuring that $\tau<\epsilon$, the second step \eqref{opsplit2} can be discretized implicitly
\begin{multline*}
  (id-\epsilon\partial_x^2)\frac{A-\hat A^{k+1}}{\tau} = -\partial_x^2 A\\
   \Rightarrow
 A^{k+1} = \left(id - (\epsilon-\tau)\partial_{xx}^2\right)^{-1}(id-\epsilon\partial_{xx}^2)\hat A^{k+1} \,.
\end{multline*}
The operator on the right-hand side is positive definite (for $\tau<\epsilon$) and is straightforward to calculate numerically with spectral methods (FFT).
Moreover, all the eigenvalues of the operator have modulus smaller than 1, which reaffirms the control of high frequency modes by the regularization.

\begin{figure}[h]
  \centering
  \includegraphics[width=0.48\textwidth]{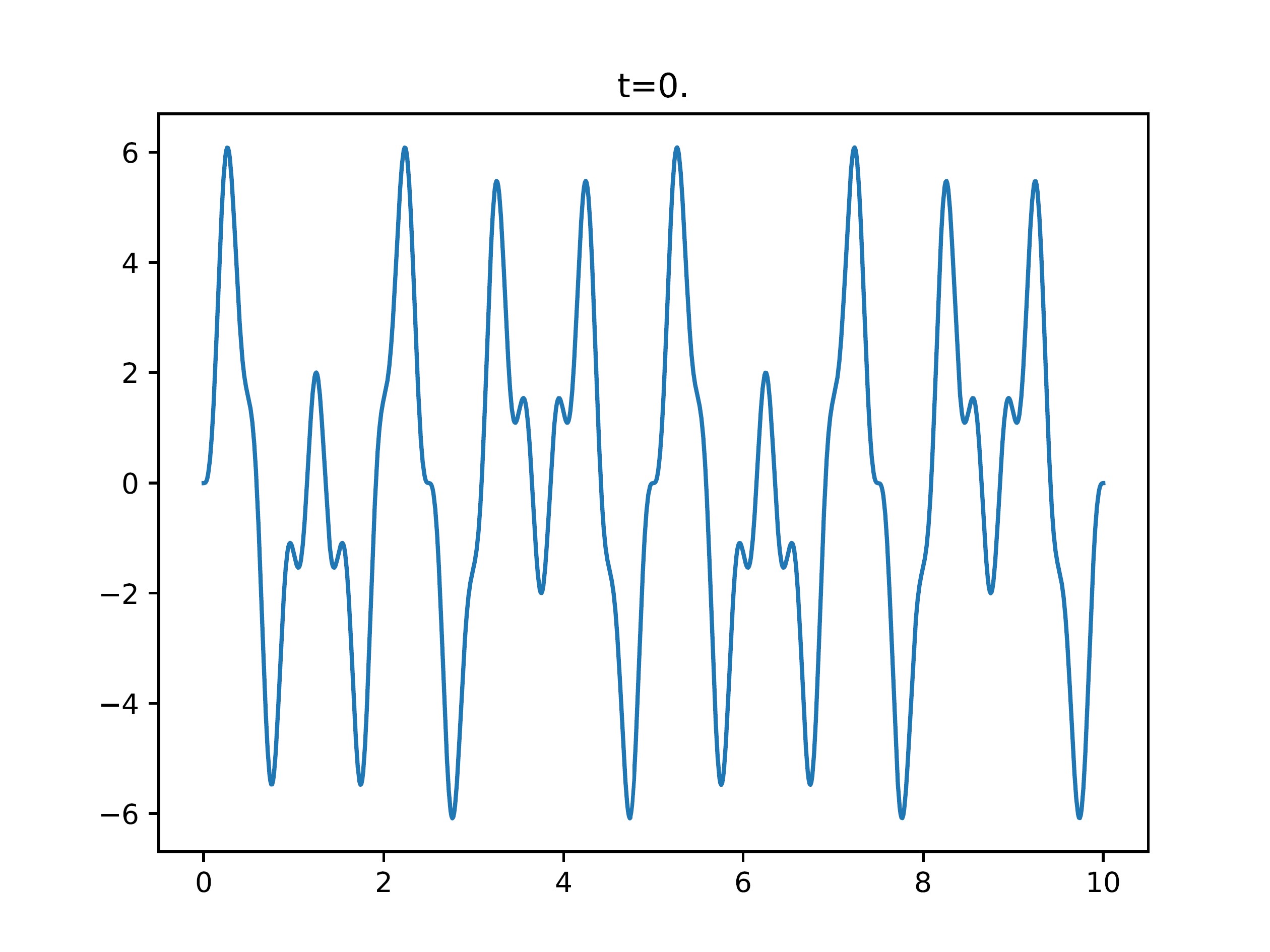} \includegraphics[width=0.48\textwidth]{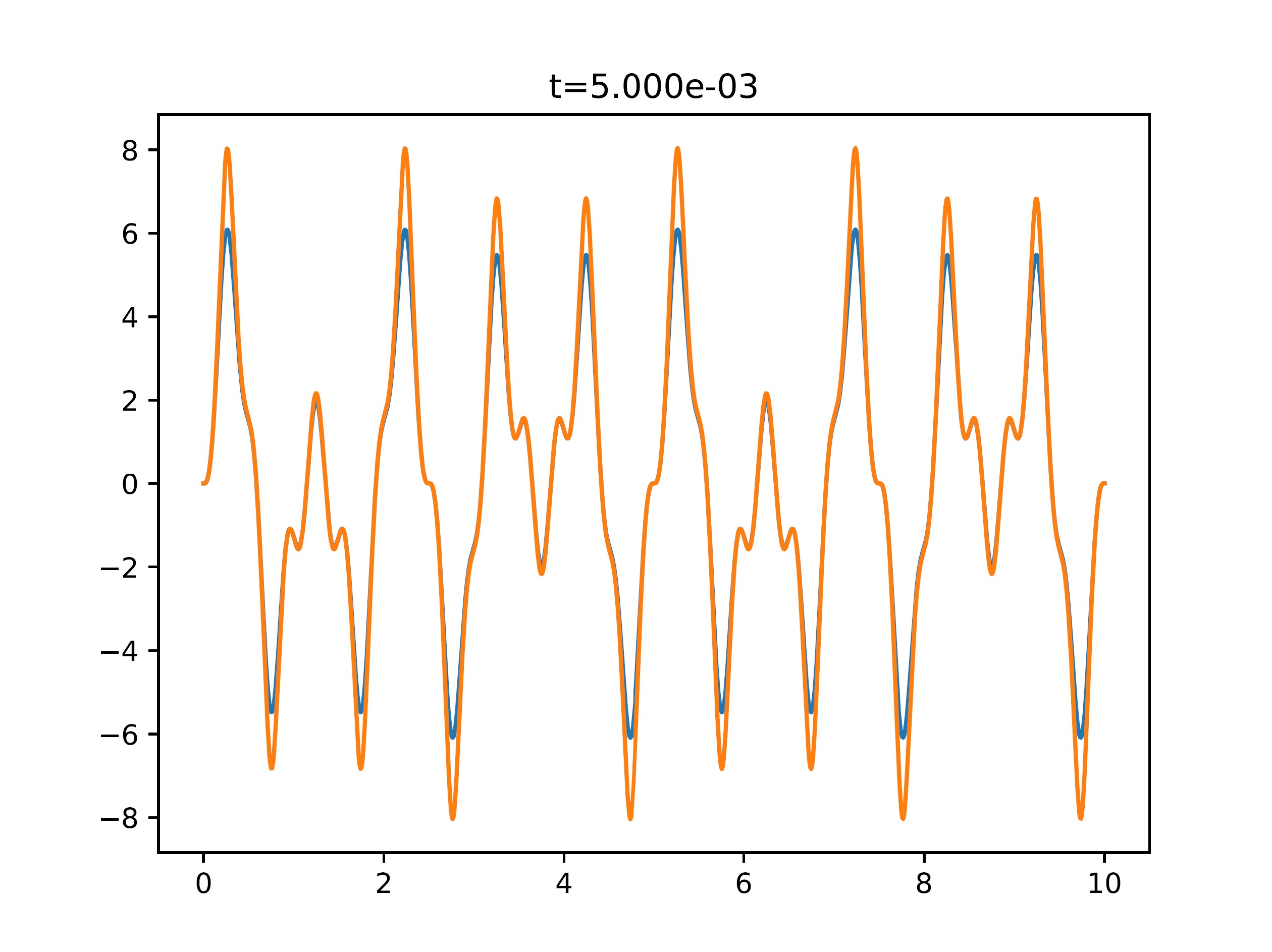}\\
  \includegraphics[width=0.48\textwidth]{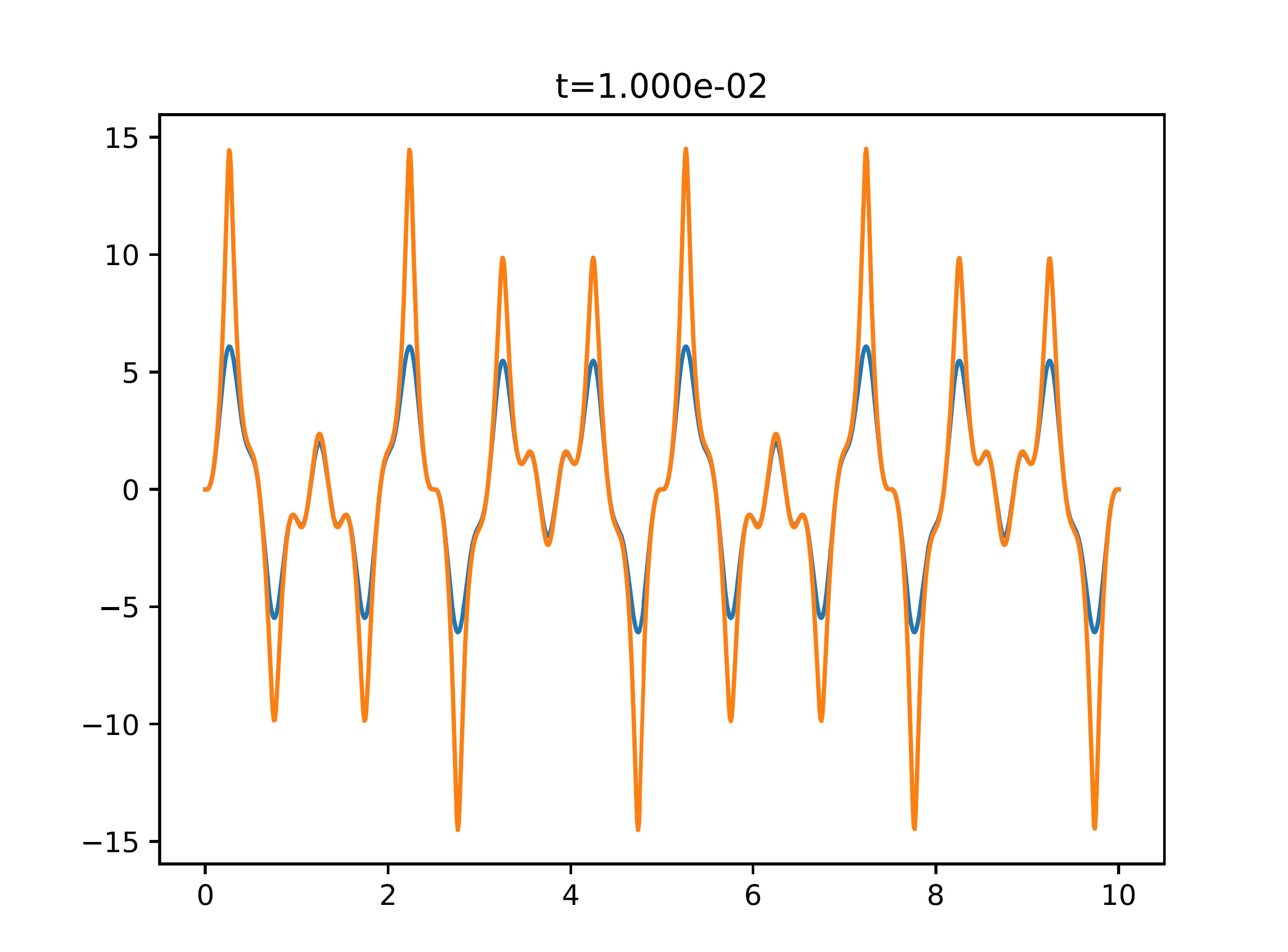} \includegraphics[width=0.48\textwidth]{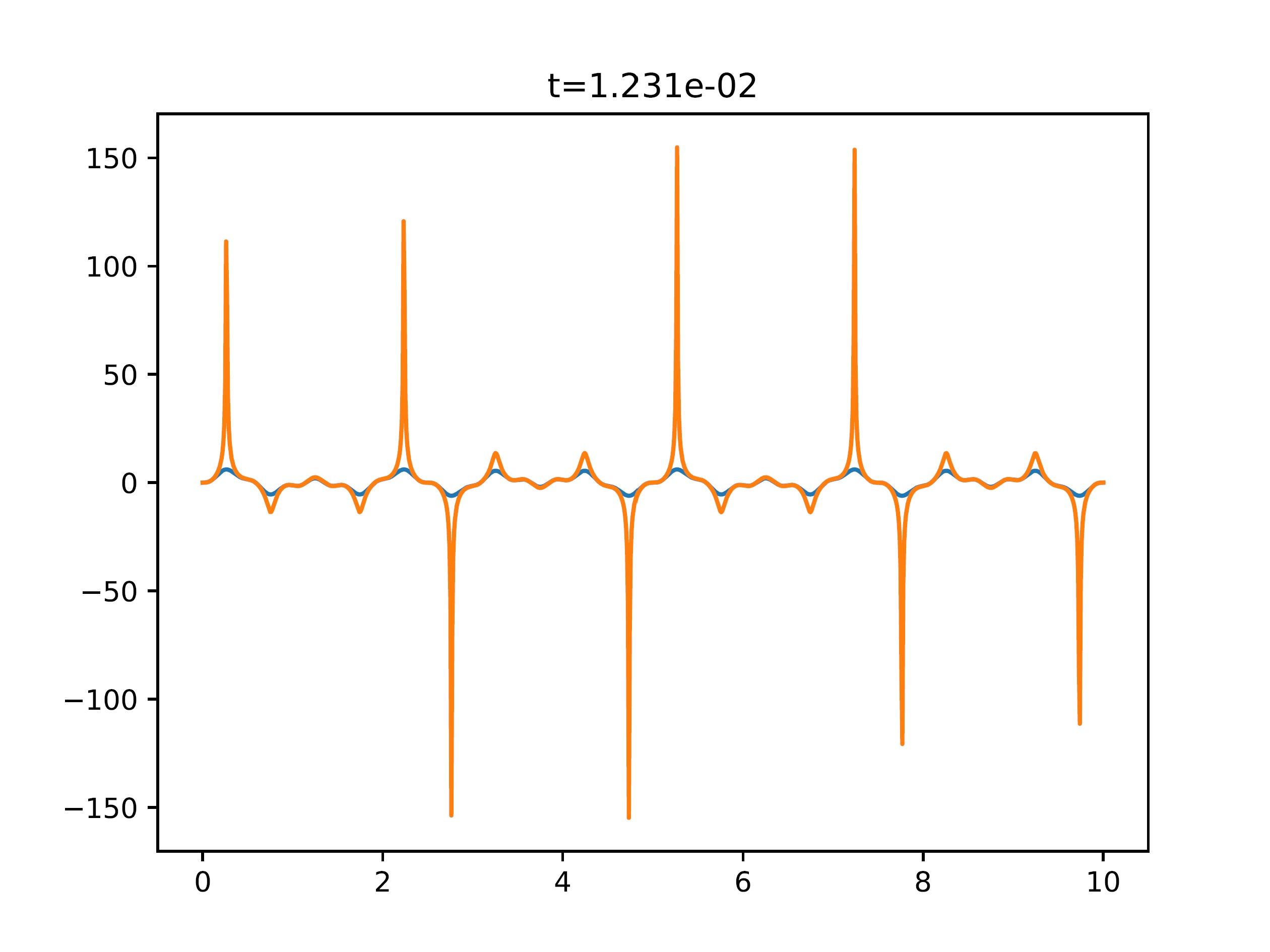}\\
  \includegraphics[width=0.48\textwidth]{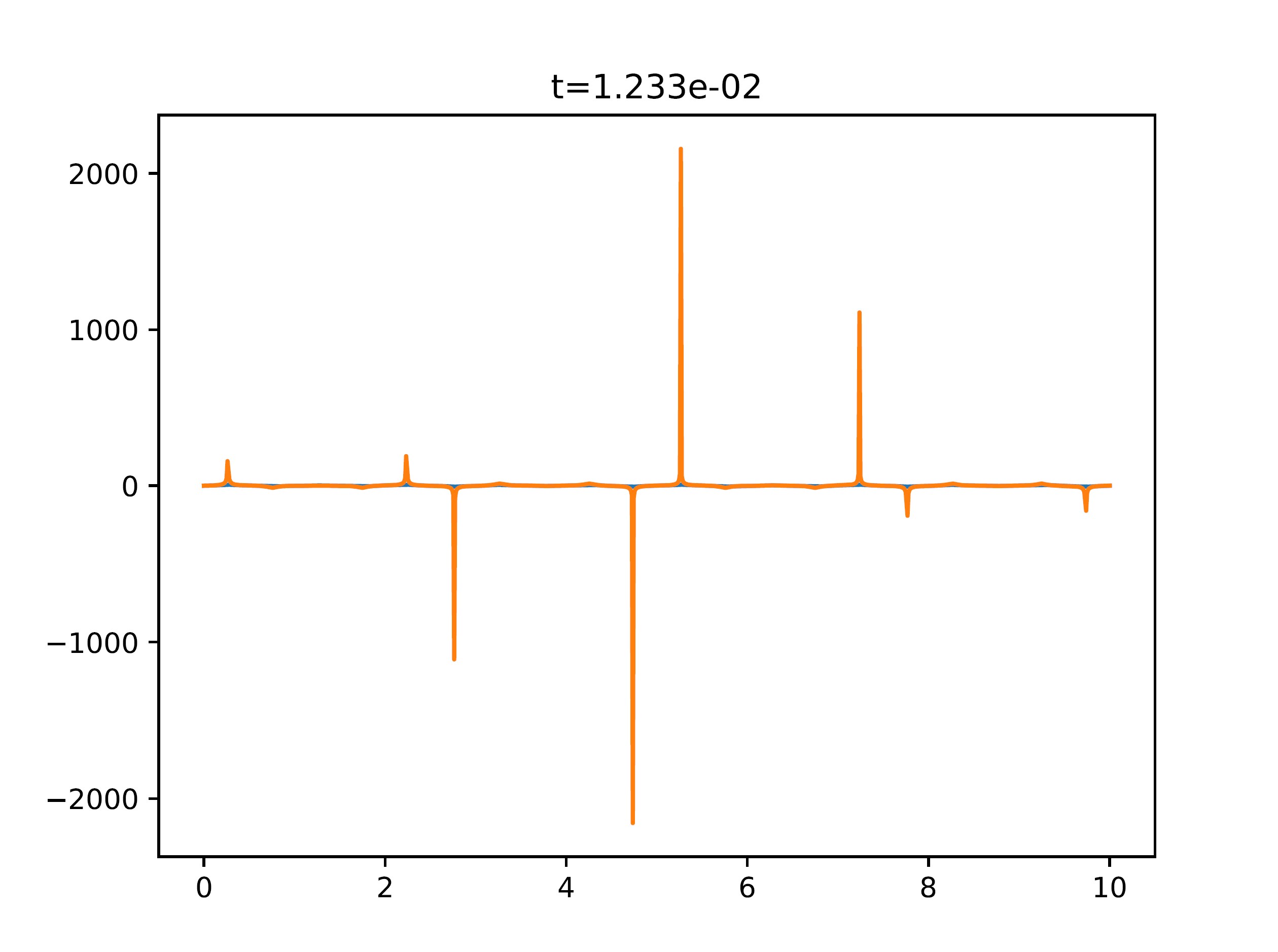} \includegraphics[width=0.48\textwidth]{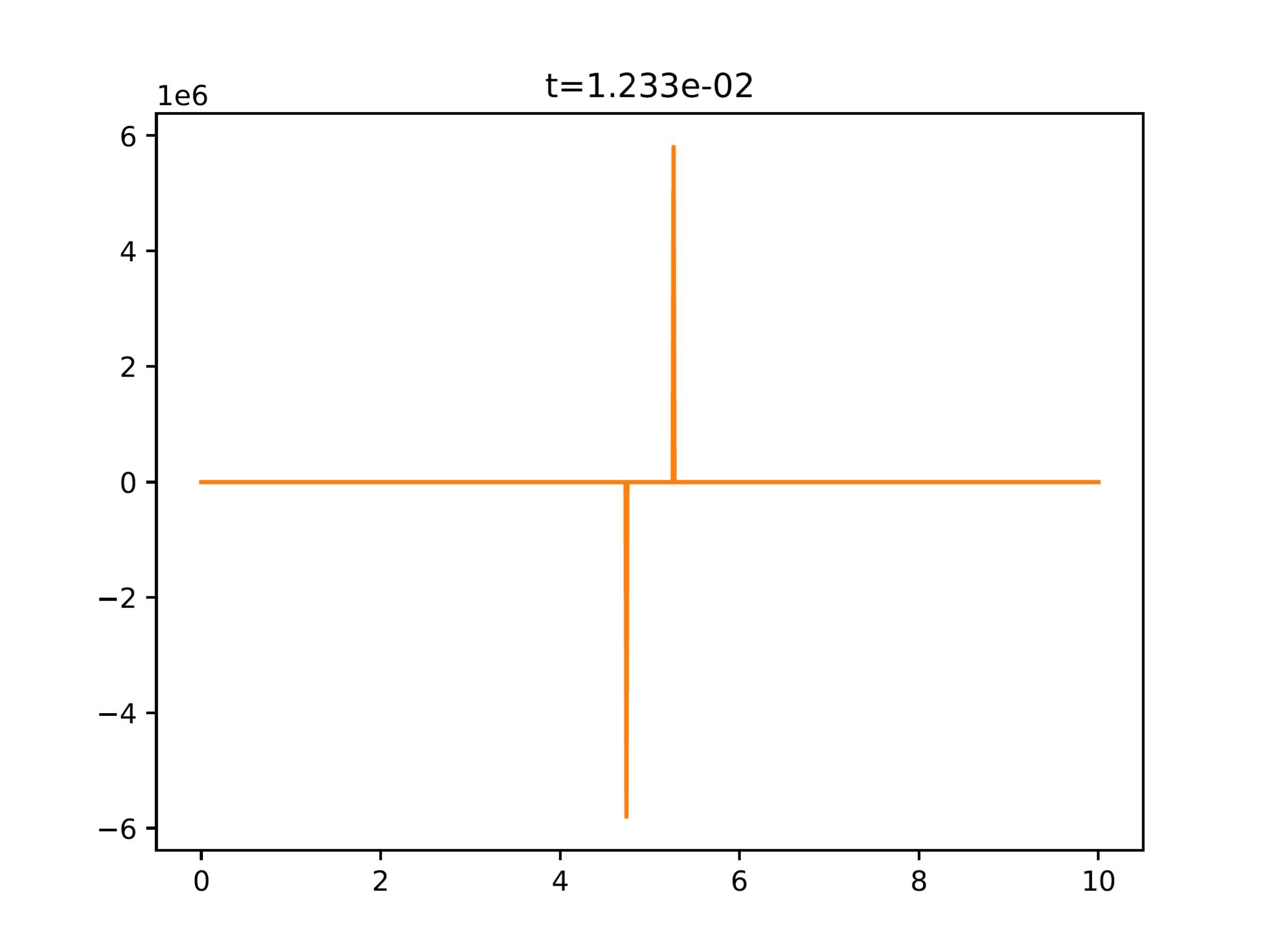}\\
  \caption{\footnotesize{Evolution of an oscillating initial condition under the regularized PDE $\eqref{Aeps}$ with $\epsilon=0.1$, 
  via the operator splitting scheme \eqref{opsplit1}-\eqref{opsplit2}. The evolution is tracked up to the time $t_{blow}$, 
  when the adaptive time step $\tau\approx 0$ up to machine precision ($\sim 10^{-16}$).}}
  \label{fig:blowup}
 \end{figure}

\begin{figure}[h]
  \centering
  \hspace{10pt}\includegraphics[height=150pt]{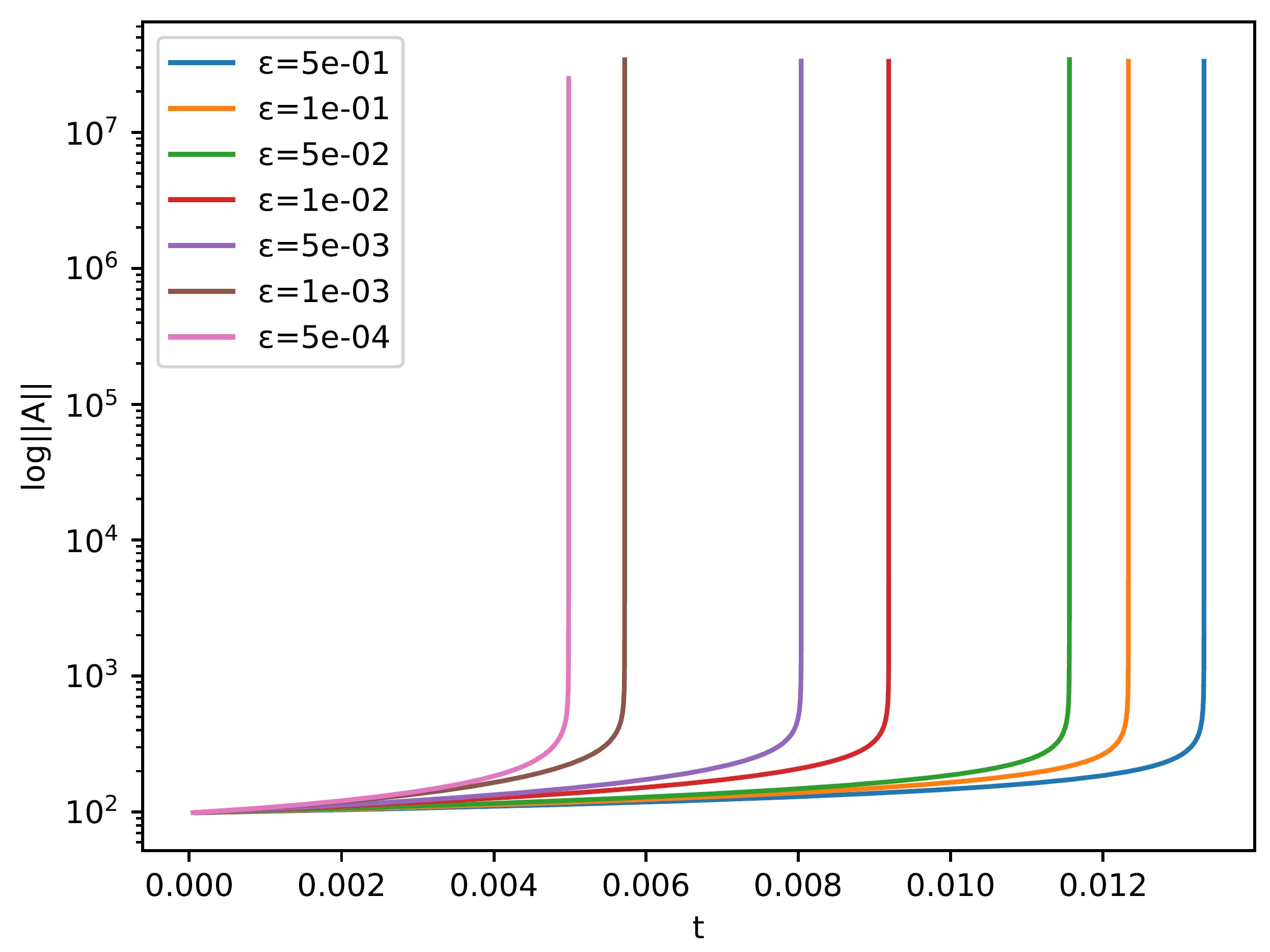}\\  \includegraphics[height=155pt]{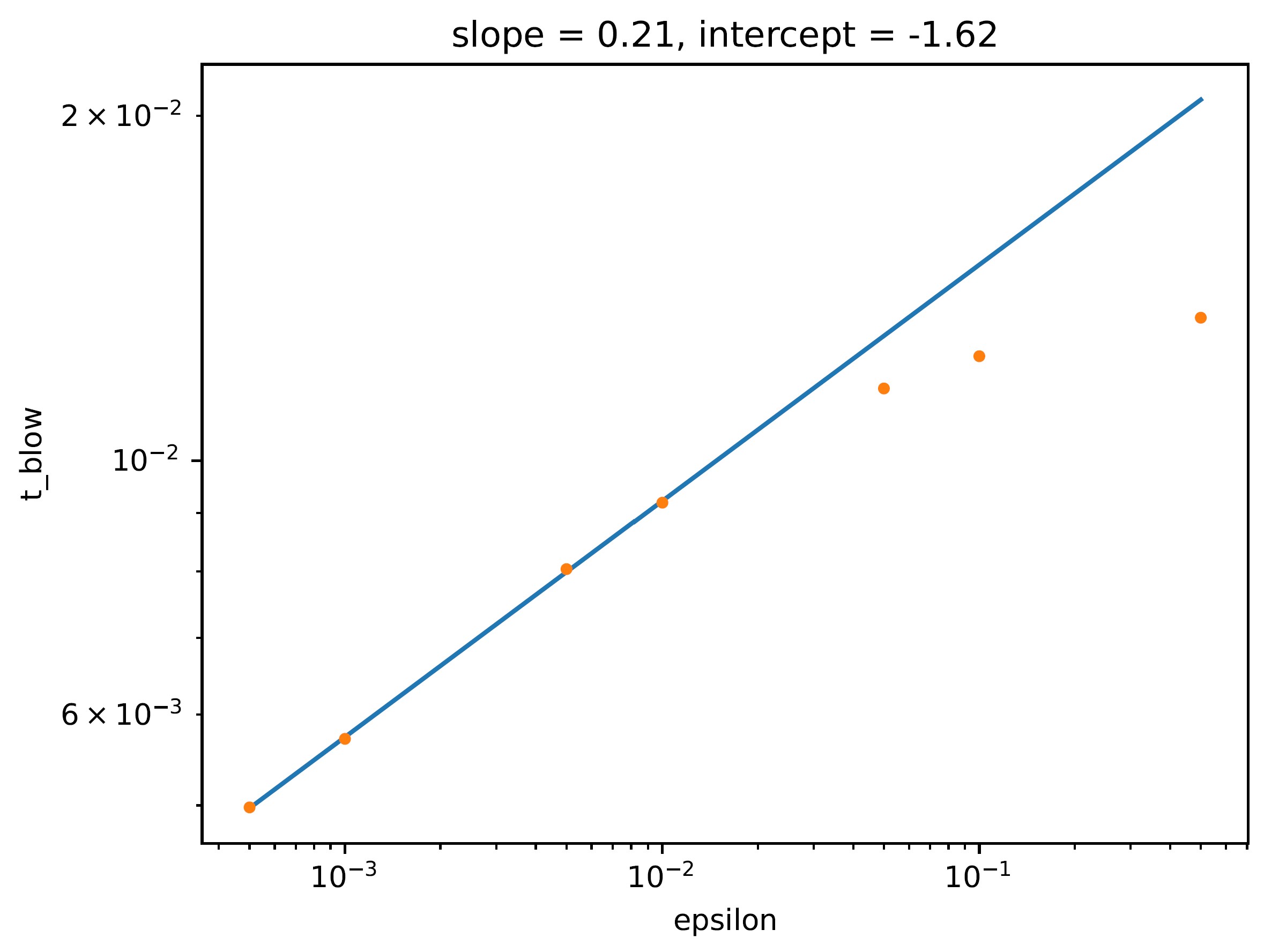}\\
  \caption{\footnotesize{Effect of varying the regularizing parameter $\epsilon$. 
  In the top plot, we see the norm $\lVert A\rVert$ blowing up at progressively earlier times for smaller values of $\epsilon$.
  In the bottom plot, we look at the log-log plot of blow-up time vs the parameter $\epsilon$, which appears to be $t_{blow} \sim \epsilon^{0.21}$ as $\epsilon\rightarrow 0$.
  }}
  \label{fig:blowup_analysis}
 \end{figure}

In figure \ref{fig:blowup}, we show the results of applying the scheme to the initial condition $A_0(x) = 5\sin^3(20\pi x/L) + 2\sin^3(12\pi x/L) - \sin^3(4\pi x/L)$
with $L=10$. We discretize the interval $[0,L]$ with 1000 uniformly distributed points, and we track the evolution of $A$ up to the time $t_{blow}$, 
when the adaptive time step $\tau\approx 0$ up to machine precision ($\sim 10^{-16}$). 
We repeat the numerical experiment with different values for the regularizing parameter $\epsilon\in\{0.5,0.1,0.05,0.01,0.005,0.001,0.0005\}$, with figure 2 corresponding to $\epsilon=0.1$.
The results are presented in figure \ref{fig:blowup_analysis}, where we see that the numerical solution blows up at earlier times $t_{blow}$  for smaller values of $\epsilon$,
with a measured rate of $t_{blow}\sim \epsilon^{0.21}$ as $\epsilon\rightarrow 0$.

\newpage
\section{Determining functionals}\label{sec:functionals}
In this section, we will show that finitely many linear functionals defined on the phase space uniquely determine the asymptotic behavior of solutions.\\
To this end, first we derive a uniform a priory estimates of solutions to the problem in the phase space $V^1:=H^1_0(\Om)\times H^1_0(\Om)$.\\
Multiplying  the equation \eqref{a0} with $-\Delta A^{\ast}$ in $L^2(\Om)$ and employing the Cauchy-Schwarz inequality we get
\begin{multline}\label{a0e}
\frac{\tau}2 \frac d{dt} \|\Nx A\|^2+\|\Dx A\|^2+(\phi^{2},\lvert\nabla A\rvert^{2})+(|A|^2,|\Nx A|^2)\\
\le \|\Nx A\|^2+2\lvert( A\partial_{x}A^{\ast},\phi\partial_{x}\phi)\rvert+2\lvert(A\partial_{y}A^{\ast},\phi\partial_{y}\phi)\rvert\\
+2c_1\|\phi \partial_{y}A\|\|\Dx A\|+\lvert\beta\rvert\lVert A\partial_{y}\phi\rVert\lVert\Delta A\rVert
\end{multline}
Utilizing the inequality
\be\label{intin}
\|\Nx v \|^2\le \|v\|\|\Dx v\|, \ \ \forall v\in H^2(\Om)\cap H_0^1(\Om),
\ee
the Ladyzhenskaya inequality
\be\label{Lad4}
\|v\| _{L^4(\Om)}^2\le \sqrt{2}\|v\|\|\Nx v\|, \ \ \forall v\in H_0^1(\Om),
\ee
the interpolation inequality,
\begin{equation}\label{ineq:interpolation}
\lVert\nabla v\rVert_{L^{4}(\Omega)}\leq d_{0}\lVert\Delta v\rVert^{\frac{1}{2}}\lVert \nabla v\rVert^{\frac{1}{2}},\quad \forall v\in H^{2}(\Omega)\cap H^{1}_{0}(\Omega),
\end{equation}
the Young's inequality
\be\label{Young}
ab\leq \frac{\eb}p a^p+\frac1{q\eb^{\frac1{p-1}}}b^q, \mbox{for all} \
a,b,\eb>0,\ \mbox{with}\ q=p/(p-1), 1<p<\infty,
\ee
we derive the following inequality from \eqref{a0e}
\begin{multline}\label{sec5:ineq0}
\frac{\tau}2 \frac d{dt} \|\Nx A\|^2+\|\Dx A\|^2+(\phi^{2},\lvert\nabla A\rvert^{2})+(|A|^2,|\Nx A|^2)\\ \leq3\lVert\Delta A\rVert^{2}+\frac{1}{4\varepsilon}\lVert A\rVert^{2}+2\vert(A\partial_{x}A^{\ast},\phi\partial_{x}\phi)\rvert\\
+2\lvert(A\partial_{y}A^{\ast},\phi\partial_{y}\phi)\rvert+\frac{c_{1}^{2}}{\varepsilon}\lVert\phi\partial_{y}A\rVert^{2}+\frac{\beta^{2}}{4\varepsilon}\lVert A\partial_{y}\phi\rVert^{2}.
\end{multline}
We obtain the following estimates for the terms on the right hand side of \eqref{sec5:ineq0} by using the following inequality
\be\label{Ellip}
\nu_0 \|\Dx \phi\|^2\le \|\mathbb L \phi\|^2\le \nu_1 \|\Dx \phi\|^2,\quad \nu_0>0,\quad \nu_1>0,
\ee
in addition to \eqref{intin}, \eqref{Lad4}, \eqref{ineq:interpolation} and \eqref{Young}:
\begin{multline*}
2\lvert(A\partial_{x}A^{\ast},\phi\partial_{x}\phi)\rvert+2\lvert(A\partial_{y}A^{\ast},\phi\partial_{y}\phi)\rvert\\
\leq 4\lVert A\rVert_{L^{4}(\Omega)}\lVert \phi\rVert_{L^{4}(\Omega)}\lVert\nabla \phi\rVert_{L^{4}(\Omega)}\lVert\nabla A\rVert_{L^{4}(\Omega)}\quad\quad\quad\quad\quad\quad\\
\leq 8 \lVert A\rVert^{\frac{1}{2}}\lVert\nabla A\rVert^{\frac{1}{2}}\lVert\phi\rVert^{\frac{1}{2}}\lVert\nabla\phi\rVert^{\frac{1}{2}}\lVert\nabla\phi\rVert^{\frac{1}{2}}\lVert\Delta\phi\rVert^{\frac{1}{2}}\lVert\nabla A\rVert^{\frac{1}{2}}\lVert\Delta A\rVert^{\frac{1}{2}}\\
\leq\frac{\varepsilon}{4}\lVert\Delta\phi\rVert^{2}+\frac{\varepsilon}{4}\lVert\Delta A\rVert^{2}+C(\varepsilon)\lVert A\rVert^{2}\lVert\phi\rVert^{2}\lVert\nabla A\rVert^{4}+\frac{\lVert\nabla \phi\rVert^{4}}{4},\quad\quad\quad\
\end{multline*}
\begin{multline*}
\frac{c_{1}^{2}}{\varepsilon}\lVert\phi\partial_{y}A\rVert^{2}\leq\frac{c_{1}^{2}}{\varepsilon}\lVert\phi\rVert^{2}_{L^{4}(\Omega)}\lVert\nabla A\rVert^{2}_{L^{4}(\Omega)}\\
\leq\frac{c_{1}^{2}d_{0}^{2}}{\varepsilon}\lVert\phi\rVert^{2}_{L^{4}(\Omega)}\lVert\Delta A\rVert\lVert\nabla A\rVert\leq\frac{c_{1}^{2}d_{0}^{2}}{\varepsilon}\lVert\phi\rVert^{2}_{L^{4}(\Omega)}\lVert A\rVert^{1/2}\lVert\Delta A\rVert^{3/2}\\
\leq \frac{3\varepsilon}{4}\lVert\Delta A\rVert^{2}+\frac{c_{1}^{8}d_{0}^{8}}{4\varepsilon^{7}}\lVert\phi\rVert^{8}_{L^{4}(\Omega)}\lVert A\rVert^{2}\\
\leq\frac{3\varepsilon}{4}\lVert\Delta A\rVert^{2}+C(\varepsilon, c_{1},d_{0})\lVert\phi\rVert^{4}\lVert\nabla\phi\rVert^{4}\lVert A\rVert^{2}.
\end{multline*}
By following the same sequence of steps, we can similarly get
\begin{eqnarray*}
\frac{\beta^{2}}{4\varepsilon}\lVert A\partial_{y}\phi\rVert^{2}&\leq&\frac{\beta^{2}}{4\varepsilon}\lVert A\rVert^{2}_{L^{4}(\Omega)}\lVert\nabla \phi\rVert^{2}_{L^{4}(\Omega)}\\
&\leq&\frac{3\varepsilon}{4}\lVert\Delta\phi\rVert^{2}+C(\varepsilon,\beta,d_{0})\lVert A\rVert^{4}\lVert\nabla A\rVert^{4}\lVert\phi\rVert^{2}.
\end{eqnarray*}
Thus we infer from \eqref{sec5:ineq0} that
\begin{multline}\label{sec5:ineq1}
\frac{\tau}2 \frac d{dt} \|\Nx A\|^2+(1-4\varepsilon)\|\Dx A\|^2+(\phi^{2},\lvert\nabla A\rvert^{2})+(|A|^2,|\Nx A|^2)\\ \leq\frac{\varepsilon}{4\nu_{0}}\lVert\mathbb{L}\phi\rVert^{2}+\frac{1}{4\varepsilon}\lVert A\rVert^{2}+C(\varepsilon)\lVert A\rVert^{2}\lVert\phi\rVert^{2}\lVert\nabla A\rVert^{4}+\frac{\lVert\nabla\phi\rVert^{4}}{4}\\
+C(\varepsilon,c_{1},d_{0})\lVert\phi\rVert^{4}\lVert\nabla\phi\rVert^{4}\lVert A\rVert^{2}+C(\varepsilon,\beta,d_{0})\lVert A\rVert^{4}\lVert\nabla A\rVert^{4}\lVert\phi\rVert^{2}.
\end{multline}

\noindent Next, we multiply the equation \eqref{a1} in $L^2(\Om)$ with $-\mathbb{L}\phi$:
\begin{multline}\label{a0e3}
\frac12 \frac d{dt} (\mathbb L \phi,\phi)+\|\mathbb L \phi\|^2 +h(\mathbb L \phi,\phi)\\
\le (|\phi||A|^2,|\mathbb L \phi|)+|c_2|(|A|\lvert\partial_{y}A\rvert, |\mathbb L \phi|)
\end{multline}
Arguing similarly to the calculation of the estimate \eqref{sec5:ineq1} and employing the interpolation inequality
$$
\|v\|_{L^6(\Om)}\le d_{1} \|\Nx v\|^{\frac{2}{3}}\| v\|^{\frac{1}{3}}, 
$$ 
we obtain:
\begin{eqnarray*}
(\lvert\phi\rvert\lvert A\rvert^{2},\lvert\mathbb{L}\phi\rvert)&\leq&\varepsilon\lVert\mathbb{L}\phi\rVert^{2}+\frac{1}{4\varepsilon}\lVert\phi\rVert_{L^{6}}^{2}\lVert A\rVert_{L^{6}}^{4}\\
&\leq&\varepsilon\lVert\mathbb{L}\phi\rVert^{2}+\frac{d_{1}^{6}}{4\varepsilon}\lVert\nabla\phi\rVert^{4/3}\lVert\phi\rVert^{2/3}\lVert\nabla A\rVert^{8/3}\lVert A\rVert^{4/3}\\
&\leq&\varepsilon\lVert\mathbb{L}\phi\rVert^{2}+C(\varepsilon,d_{1})(\lVert\phi\rVert^{2}\lVert\nabla\phi\rVert^{4}+\lVert A\rVert^{2}\lVert\nabla A\rVert^{4}),
\end{eqnarray*}
\begin{eqnarray*}
\lvert c_{2}\rvert(\lvert A\rvert\lvert\partial_{y}A\rvert,\lvert \mathbb{L}\phi\rvert)&\leq& \varepsilon\lVert\mathbb{L}\phi\rVert^{2}+\frac{c_{2}^{2}}{4\varepsilon}(\lvert A\rvert^{2},\lvert\partial_{y}A\rvert^{2})\\
&\leq&\varepsilon\lVert\mathbb{L}\phi\rVert^{2}+\frac{c_{2}^{2}}{4\varepsilon}\lVert A\rVert_{L^{4}(\Omega)}\lVert\partial_{y}A\rVert^{2}_{L^{4}(\Omega)}\\
&\leq&\varepsilon\lVert\mathbb{L}\phi\rVert^{2}+\frac{\sqrt{2}c_{2}^{2}d_{0}^{2}}{4\varepsilon}\lVert A\rVert\lVert\nabla A\rVert^{2}\lVert\Delta A\rVert\\
&\leq&\varepsilon(\lVert\mathbb{L}\phi\rVert^{2}+\lVert\Delta A\rVert^{2})+C(\varepsilon,c_{2},d_{0})\lVert A\rVert^{2}\lVert\nabla A\rVert^{2}.
\end{eqnarray*}
By using the above estimates, we infer from \eqref{a0e3}
\begin{multline}\label{afi}
\frac{1}{2}\frac{d}{dt}(\mathbb{L}\phi,\phi)+(1-2\varepsilon)\lVert\mathbb{L}\phi\rVert^{2}-\varepsilon\lVert\Delta A\rVert^{2}+h(\mathbb{L}\phi,\phi)\\
\leq C(\varepsilon,d_{1})(\lVert\phi\rVert^{2}\lVert\nabla\phi\rVert^{4}+\lVert A\rVert^{2}\lVert\nabla A\rVert^{4})+C(\varepsilon, c_{2},d_{0})\lVert A\rVert^{2}\lVert\nabla A\rVert^{4}.
\end{multline}
Finally adding \eqref{sec5:ineq1} and \eqref{afi} we get
\begin{multline}\label{sec5:ineq2}
\frac d{dt}\left[\frac{\tau}2 \|\Nx A\|^2+\frac12  (\mathbb L \phi,\phi) \right]
+(1-5\varepsilon)\|\Dx A\|^2\\
+(1-2\varepsilon-\frac{\varepsilon}{4\nu_{0}})\lVert\mathbb{L}\phi\rVert^{2}+(\phi^{2},\lvert\nabla A\rvert^{2})+(\lvert A\rvert^{2},\lvert\nabla A\rvert^{2})\\
\leq\frac{1}{4\varepsilon}\lVert A\rVert^{2}+C(\varepsilon)\lVert A\rVert^{2}\lVert \phi\rVert^{2}\lVert\nabla A\rVert^{4}+\frac{\lVert\nabla\phi\rVert^{4}}{4}+C(\varepsilon,c_{1},d_{0})\lVert\phi\rVert^{4}\lVert\nabla\phi\rVert^{4}\lVert A\rVert^{2}\\
+C(\varepsilon,\beta,d_{0})\lVert A\rVert^{4}\lVert\nabla A\rVert^{4}\lVert\phi\rVert^{2}+C(\varepsilon,d_{1})(\lVert\phi\rVert^{2}\lVert\nabla\phi\rVert^{4}+\lVert A\rVert^{2}\lVert\nabla A\rVert^{4})\\
+C(\varepsilon,c_{2},d_{0})\lVert A\rVert^{2}\lVert\nabla A\rVert^{4}
\end{multline}
Since the semigroup generated by the problem is dissipative in $V^{0}$ there exist $T_{0}>0$ and $M_{0}>0$ such that
\begin{equation*}
\lVert A(t)\rVert^{2}+\lVert\phi(t)\rVert^{2}\leq M_{0},\quad \forall t\geq T_{0}.
\end{equation*}
Thus, by choosing $\varepsilon=\min\{\frac{1}{10},\frac{2\nu_{0}}{8\nu_{0}+1}\}$ we obtain from \eqref{sec5:ineq2} that
\begin{multline*}
\frac{d}{dt}[\tau\lVert\nabla A\rVert^{2}+(\mathbb{L}\phi,\phi)]+\lVert\Delta A\rVert^{2}\lVert\mathbb{L}\phi\rVert^{2}\\
\leq\frac{1}{4\varepsilon}\lVert A\rVert^{2}+C(\varepsilon, M_{0},\beta,d_{0},d_{1})(\lVert\nabla A\rVert^{4}+\lVert\nabla\phi\rVert^{4}).
\end{multline*}
Finally, by using the uniform Gronwall Lemma \ref{uGr}, we infer from the last inequality that there exists 
$R_0>0$,  depending only on the parameters of the system, such that
\be\label{phia}
\|\Nx A\|^2+ \|\Nx \phi\|^2\le R_0, \ \ \forall  t\ge T_1,
\ee
i.e., the semigroup generated by the  problem is bounded dissipative in $V^1=H_0^1(\Om)\times H_0^1(\Om)$.\\
It remains to show that the condition \eqref{Lip1} is satisfied, i.e.,
\be\label{fja}
\|f_j(a,\Phi)\|_{H^{-\gamma}} \le M(R_0)(\|\Nx a\|+\|\Nx \Phi\|), \ \ j=1,2,
\ee
with some $\gamma \in (0,1)$.

In fact, employing the estimate \eqref{phia}, the following  inequality  (see, e.g., \cite{BV},  Chp. 1)
$$
\|v\|_{H^{-\delta}(\Om)}\le  d_0 \|v\|_{L^{2/(1+\delta)}(\Om)}
$$
which is valid for $\delta \in (0,1)$, and the fact that $H^1(\Om)$  is continuously embedded into $L^{\frac2\delta}(\Om)$ we can bound the last term on the right hand side of \eqref{f11}
\begin{multline*}
\|i\beta[\t A\p_y\Phi]\|_{H^{-\delta}(\Om)}\le  d_0 |\beta|\|\t A\p_y\Phi\|_{L^{2/(1+\delta)}}\\
\le d_0 |\beta|\|\t A\|_{L^{2/\delta}}\lVert\partial_{y}\Phi\rVert\le  Cd_0\|\Nx \t A\|
\|\p_y\Phi\|\le C(R_0)\|\Nx \Phi\|.
\end{multline*}
Similarly we can show that 
$$
\|i\beta[a\p_y\phi]\|_{H^{-\delta}(\Om)}\le C(R_0)\|\Nx a\|.
$$
Utilizing the same arguments we  can get necessary bounds  for all terms on the right-hand side of \eqref{f11} and \eqref{f22} and obtain the desired estimates \eqref{fja}.
Thus  thanks to  the Theorem \ref{Det1} we have

\begin{theorem} \label{DF2}The  set $\LL $ of functionals on $V^{1}$ 
is asymptotically  $(V^1, V^{0},\mathcal{W})$  determining functionals for the problem \eqref{a0}-\eqref{a3} provided the defect $\eb(V^{1}, V^{0})$ is small enough.
\end{theorem} 

\begin{remark} Similar to the proof of the Theorem \ref{DF2} we can prove that the  set $\LL $ of functionals on $V^{1}$ 
is asymptotically  $(V^1, V^{0}, \mathcal{W})-$  determining functionals for the problem 1D problem
\be\label{1Df}
\begin{cases}
 \tau \partial_{t}A=A+ \p_x^2A - \phi^2A-|A|^2A, \ x\in (0,L), t>0,\\
\partial_{t}\phi=\p_x^2\phi-h\phi +\phi\lvert A\rvert^{2}, \ x\in (0,L), t>0,\\
 A\big|_{x=0}= A\big|_{x=L}=\phi\big|_{x=0}=\phi\big|_{x=L}=0,\\ 
 A\big|_{t=0}=A_0,\quad \phi\big|_{t=0}=\phi_0
\end{cases}
\ee
provided the defect $\eb(V^{1} ,{\bf H})$ is small enough. Let us note that dissipativity of the semigroup $S(t),t\ge0$ generated by the problem \eqref{1Df} , as well as existence of 
and exponential attractor for the semigroup $S(t),t\ge0$ is established in \cite{KKV2021}.

\end{remark}

\end{document}